\documentclass{article}
\usepackage{amsmath}  
\usepackage{amssymb}

\newcommand{\sn}{\smallskip\noindent}
\newcommand{\mn}{\medskip\noindent}
\newcommand{\bn}{\bigskip\noindent}

\newcommand{\qed}{\vrule height 1.2ex width 1.1ex depth -.1ex}

\newcommand{\A}{\mathcal{A}}

\newcommand{\R}{\mathcal{R}}
\newcommand{\C}{\mathcal{C}}

\newcommand{\Oo}{\mathcal{O}}

\newcommand{\id}{\mathrm{id}}
\newcommand{\Lin}{\mathrm{Lin}}

\newcommand{\bbbz}{\mathbb{Z}}
\newcommand{\bbbr}{\mathbb{R}}
\newcommand{\bbbc}{\mathbb{C}}
\newcommand{\br}{{\bf r}}
\newcommand{\bs}{{\bf s}}
\newcommand{\bc}{{\bf c}}

\newcommand{\BWM}{\mathrm{BWM}}
 
\newcommand{\Mor}{\mathrm{Mor}}

\author{\sc Konrad Schm\"udgen}
\title{{\bf ON COQUASITRIANGULAR BIALGEBRAS}}
\date{\small \sc Universit\"at Leipzig,
Fakult\"at f\"ur Mathematik und Informatik,
Augustusplatz 10/11, D-04109 Leipzig, Germany\\
E-mail: schmuedg@mathematik.uni-leipzig.de}

\begin{document}

\maketitle

\bn
{\bf 1. Introduction} 

\mn
Coquasitriangularity is one of the most fundamental concepts in quantum group theory. Some early papers on this notion are \cite{Lawson}, \cite{Schauenburg}, \cite{Hayashi}, \cite{Doi}, \cite{Majid}. A bialgebra or a Hopf algebra is called {\it coquasitriangular} if it is equipped with a universal $r$-form. The latter concept is the main object of study in this paper.

\bn
{\bf Definition 1.1} Let $\A$ be a bialgebra. A {\it universal r-form} on $\A$ is a linear functional on $\A\otimes \A$ which is invertible with respect to the convolution multiplication and satisfies the following three conditions for arbitrary $a,b,c\in\A$:

\sn
(CQT.1)\quad $\br (c\otimes ab)=\br (c_{(1)} \otimes b) \br(c_{(2)}\otimes a),$\\
(CQT.2)\quad $\br(ab\otimes c)=\br(a\otimes c_{(1)}) \br(b\otimes c_{(2)})$,\\
(CQT.3)\quad $\br (a_{(1)}\otimes b_{(1)}) a_{(2)} b_{(2)}=\br(a_{(2)}\otimes b_{(2)}) b_{(1)} a_{(1)}$.

\mn
The present paper deals with three topics on coquasitriangular bialgebras. 
In Section 2 we give a characterization of a universal $r$-form in terms 
of (certain) Yetter-Drinfeld modules. In Section 3 we study the uniqueness 
of universal $r$-forms for the coordinate Hopf algebras of the quantum 
groups $GL_q(N)$, $SL_q(N)$, $O_q(N)$ and $Sp_q(N)$. Section 4 is concerned 
with some linear functionals on a coquasitriangular Hopf algebra $\A$ 
with universal $r$-form $\br$ which describe the square of the 
antipode of $\A$. To be more precise, it is known that 
$S^2=\bar{f_\br}\ast \id\ast f_\br$, where $f_\br$ is the 
linear functional on $\A$ defined by 
$f_\br(a)=\br (a_{(1)}, S(a_{(2)})), a\in\A$, and $\bar{f_\br}$ is the 
convolution inverse of $f_\br$. Among others we show that 
$F_\br :=f_\br\ast f_{\bar{\br}_{21}}$ is a character on $\A$ which 
coincides with Woronowicz$^\prime$s modular character $f_{-2}$ when 
$\A$ is the coordinate Hopf $\ast$-Algebra for the standard compact 
quantum groups $U_q(N), SU_q(N),O_q(N;\bbbr)$, and $USp_q(N)$. 
Dual versions for quasitriangular Hopf algebras of some results in 
Section 4 have been proved by V.G. Drinfeld \cite{Drinfeld}.

\mn
Let us fix some notation. We use the Sweedler notations 
$\Delta(a)=a_{(1)}\otimes a_{(2)}$ for the comultiplication 
$\Delta$ and $\delta (m)=m_{(0)}\otimes m_{(1)}$ for a right 
coaction $\delta$. The multiplication map of an algebra $\A$ is denoted by 
$m_\A$. We write $\ast$ for the convolution 
multiplication. The convolution inverse of a functional 
$f$ is denoted by $\bar{f}$.

\vfill\eject\noindent
\bn
{\bf 2. Universal ${\bf r}$-forms and Yetter-Drinfeld modules}

\mn
If not stated otherwise, $\A$ denotes a bialgebra in this section. First we 
recall a well-known definition (see [Y], [S2] or [Mo]).

\mn
{\bf Definition 2.1.} A {\it (right) Yetter-Drinfeld module} of $\A$ is a right $\A$-module and a right $\A$-comodule $M$ satisfying the compatibility condition
\begin{equation}\label{2.1}
m_{(0)}\triangleleft a_{(1)}\otimes m_{(1)} a_{(2)}=(m\triangleleft a_{(2)})_{(0)}\otimes a_{(1)} (m\triangleleft a_{(2)})_{(1)}
\end{equation}
for all $m\in M$ and $a\in \A$. Here the right action of $a\in\A$ on $m\in M$ is denoted by $m\triangleleft a$ and the right coaction $\delta$: $M\rightarrow M\otimes \A$ is expressed by the Sweedler notation $\delta (m)= m_{(0)}\otimes m_{(1)}, m\in M$.

\bn
For a right $\A$-comodule $M$, let $\C (M)=\Lin \{m^\prime (m_{(0)})~m_{(1)}| m\in M, m^\prime\in M^\prime\}$ denote the coefficient coalgebra of $M$ (see, for instance, \cite{Klimyk}, p.399).\\
For the following considerations we assume that $\br$ is a {\it convolution invertible} linear functional on the bialgebra $\A\otimes \A$. Recall that the convolution inverse of $\br$ is denoted by $\bar{\br}$.\\
Let $M$ be a right $\A$-comodule with right coaction 
$\delta (m)=m_{(0)}\otimes m_{(1)}$.  For $m\in M$ and $a\in\A$, we define
\begin{align*}
m\triangleleft_1 a&=\br(m_{(1)}\otimes a) m_{(0)} ~ {\rm and}~
~ m\triangleleft_2 a=\bar{\br}(a\otimes m_{(1)})m_{(0)}.
\end{align*}
Let $M_i, i=1,2$, denote the right $\A$-comodule $M$ equipped with the 
mapping $M\times \A\ni (m,a)\rightarrow m\triangleleft_i a\in M$. If 
$\br$ is a universal $r$-form, then it is known (see [Mo], Example 10.6.14) 
that $M_1$ is a Yetter-Drinfeld module. 
We shall strengthen this fact and give a characterization 
of universal $r$-forms 
in this manner.

\mn
{\bf Lemma 2.2} (i) {\it $M_1$ is a Yetter-Drinfeld module of $\A$ with 
right action $\triangleleft_1$ if and only if (CQT.1) holds for all 
$a,b\in\A$ and  $c\in\C (M)$ and (CQT.3) holds for all $a\in\C (M)$ and $b\in\A$.}\\
(ii) {\it  $M_2$ is a Yetter-Drinfeld module of $\A$ with right action 
$\triangleleft_2$ if and only if (CQT.2) is fulfilled for all $a,b \in\A$ and $c\in\C (M)$ and (CQT.3) holds for all $a\in\A$ and $b\in\C (M)$.}

\mn
{\bf Proof.} We prove the assertion for $M_1$.
It is obvious that the condition $(m\triangleleft_1 a)\triangleleft_1 b=m \triangleleft_1 ab$ is equivalent to equation (CQT.1) for $c\in\C(M)$ and $a,b\in\A$. We show that the latter implies that $m\triangleleft 1=m, m\in M$. Indeed, using the convolution inverse $\bar{\br}$ of $\br$ and condition (CQT.1) for $a=b=1$, we get
$$
\bar{\br}(c\otimes 1)=\bar{\br}(c_{(1)}\otimes 1)\br (c_{(2)}\otimes 1)\br (c_{(3)}\otimes 1)=\bar{\br} (c_{(1)}\otimes 1)\br (c_{(2)}\otimes 1)=\varepsilon (c).
$$
for $c\in\C(M)$, so that $m\triangleleft 1=\br (m_{(1)}\otimes 1)m_{(0)}=m$. The left hand and right hand sides of the Yetter-Drinfeld condition (1) are equal to 
$$
\br (m_{(1)}\otimes a_{(1)})m_{(0)}\otimes m_{(2)} a_{(2)}\qquad
{\rm and}\qquad \br (m_{(2)}\otimes a_{(2)})m_{(0)}\otimes a_{(1)}m_{(1)},
$$
respectively. Thus, (\ref{2.1}) is equivalent to (CQT.3) for $a\in \C (M)$ and $b\in\A$.\\
The assertion for $M_2$ follows by some slight modifications of the 
preceding reasoning. The relation $(m\triangleleft_2 a)\triangleleft_2 b=m\triangleleft_2 ab$, where $ a,b\in\A$ and $m\in M$, is fulfilled iff $\bar{\br} (a\otimes c_{(2)})\bar{\br}(b\otimes c_{(1)})=\bar{\br}(ab\otimes c)$ for $a,b\in\A$ and $c\in\C(M)$. The latter is obviously equivalent to the fact that $\br (a\otimes c_{(1)})\br(b\otimes c_{(2)})=\br(ab\otimes c)$ for $a,b\in\A$ and $c\in\A$. The Yetter-Drinfeld condition (\ref{2.1}) is satisfied iff $\bar{\br} (a_{(1)}\otimes b_{(1)})b_{(2)}a_{(2)}=\bar{\br} (a_{(2)}\otimes b_{(2)})a_{(1)} b_{(1)}$ for all $a\in \A$ and $b\in\C (M)$ which in turn is equivalent to equation (CQT.3) for $a\in\A$ and $b\in\C(M)$.\hfill\qed

\mn
An immediate consequence of Lemma 2.2 is 

\bn
{\bf Proposition 2.3} {\it Let $\R$ be a family of right comodules of a 
bialgebra $\A$ such that the linear span of the coefficient coalgebras 
$\C (M), M\in\R$, coincides with $\A$. Let $\br$ be a convolution 
invertible linear functional on $\A\otimes \A$. Then $\br$ is a 
universal $r$-form of $\A$ if and only if $M_1$ and $M_2$ are Yetter-Drinfeld 
modules for any comodule $M\in \R$.}

\mn
{\bf Corollary 2.4} {\it A convolution invertible linear functional 
$\br$ on $\A\otimes\A$ is a universal $r$-form of a bialgebra $\A$ if 
and only if $\A$ becomes a Yetter-Drinfeld module with the 
comultiplication as right coaction and with the mappings 
$\triangleleft_1$ and $\triangleleft_2$ as right actions.}

\mn
{\bf Remark 2.5} As shown by P. Schauenburg \cite{Schauenburg2}, the 
category of right Yetter-Drinfeld modules is equivalent to the category 
of Hopf bimodules over a Hopf algebra $\A$. Using this correspondence the 
preceding results can be reformulated in terms of Hopf bimodules of $\A$. 
In this setting they were crucial for the construction of bicovariant 
differential calculi on general coquasitriangular Hopf 
algebras (see \cite{SS}, last line on p.189, 
and  \cite{Klimyk}, Section 14.5).

\bn
{\bf 3. Uniqueness of universal $r$-forms for quantized matrix groups}

\mn
In this section let $G_q$ denote one of the quantum groups 
$GL_q(N)$, $SL_q(N)$, $O_q(N)$ or $Sp_q(N)$ and $\Oo(G_q)$ its 
coordinate Hopf algebra as defined in \cite{Faddeev}, 
\cite{Takeuchi} or \cite{Klimyk}, Chapter 9. It is known (see, for instance, \cite{Klimyk}, Theorem 10.9) that the Hopf algebra $\Oo(G_q)$ is coquasitriangular and there exists a universal $r$-form $\br_z$ of $\Oo(G_q)$ such that
\begin{equation}\label{3.1}
\br_z(u^i_j\otimes u^n_m)=zR^{in}_{jm},\quad i,j,n,m=1,{\dots},N.
\end{equation}
Here $u=(u^i_j)_{i,j=1,{\dots},N}$ is the fundamental matrix of $\Oo(G_q),R$ is the corresponding $R$-matrix as given in \cite{Faddeev}, (1.5) and (1.9), or in \cite{Klimyk}, (9.13) and (9.30), and $z$ is a fixed complex number such that $z\ne 0$ for $G_q=GL_q(N), z^N=q^{-1}$ for $G_q=SL_q(N)$ and $z^2=1$ for $G_q=O_q(N), Sp_q(N)$. 

\mn
Throughout this section we assume that $q$ is a complex number 
which is not a root of unity. (A closer look at the proof 
given below shows that it 
suffices to exclude only very few roots of unity.) We shall show 
that the above functionals $\br_z$ exhaust {\it all} universal $r$-forms of 
$\Oo(G_q)$. In order to place this result in a more general context, 
we need some preliminaries.

\bn
{\bf Definition 3.1} A {\it central bicharacter} of a bialgebra $\A$ is a convolution invertible linear functional $\bc$ on $\A\otimes\A$ such that for arbitrary $a,b,c\in\A$:\\
(CB.1)~~ $\bc(ab\otimes c)=\bc (a\otimes c_{(1)})\bc (b\otimes c_{(2)})$ and $\bc(c\otimes ab)=\bc (c_{(1)}\otimes b)\bc (c_{(2)}\otimes a),$\\
(CB.2)~~ $\bc(a\otimes b_{(1)}) b_{(2)}=\bc(a\otimes b_{(2)}) b_{(1)}$ and $\bc (a_{(1)}\otimes b) a_{(2)}=\bc(a_{(2)}\otimes b)a_{(1)}$.

\mn
Condition (CB.1) means that $\bc (\cdot \otimes\cdot)$ is a dual pairing 
of the bialgebras $\A$ and $\A^{\rm op}$, where $\A^{\rm op}$ is the 
bialgebra with the same comultiplication and the opposite 
multiplication as $\A$.\\
Condition (CB.2) is equivalent to the requirement that for any 
$a\in\A$ the linear functionals $\bc(a\otimes \cdot)$ and 
$\bc(\cdot\otimes a)$ on $\A$ are central in the dual algebra 
$\A^\prime$. A trivial example of central bicharacter is the 
counit of $\A\otimes\A$.

\mn
Suppose that $\bc$ is a central bicharacter and $\br$ is a universal 
$r$-form of $\A$. Using the condition (CQT.1)--(CQT.3) and (CB.1)--(CB.2) 
it is a straightforward matter to verify that the convolution product 
$\bc\ast\br$ is again a universal $r$-form of $\A$. 
Recall (\cite{Klimyk}, Proposition 10.2(iv)) that $\bar{\br}_{21}$ is 
another universal $r$-form of $\A$, where 
$\bar{\br}_{21}(a\otimes b):=\bar{\br} (b\otimes a),a,b\in\A$. Thus 
$\bc\ast\bar{\br}_{21}$ is also a universal $r$-form of $\A$.

\mn
{\bf Definition 3.2} A coquasitriangular bialgebra $\A$ is said to have an {\it essentially unique universal $r$-form} if there exists a universal $r$-form $\br$ of $\A$ such that for any universal $r$-form $\bs$ of $\A$ there exists a central bicharacter $\bc$ such that $\bs=\bc\ast\br$ or $\bs=\bc\ast\bar{\br}_{21}.$

\mn
{\bf Example 3.3} Let $\A=\bbbc G$ be the group Hopf algebra of an abelian group $G$. Then the universal $r$-forms of $\A$ are in one-to-one correspondence to the bicharacters of the group $G$. Since $\bbbc G$ is cocommutative and $G$ is abelian, the universal $r$-forms are precisely the central bicharacters of $\A$ and hence $\A=\bbbc G$ has obviously an essentially unique universal $r$-form.

\bn
Let us return to the coquasitriangular Hopf algebra $\Oo(G_q)$. We shall say that a complex number $\zeta$ is {\it admissible} if $\zeta\ne 0$ for $G_q=GL_q(N),\zeta^N=1$ for $G_q=SL_q(N)$ and $\zeta^2=1$ for $G_q=O_q(N), Sp_q(N)$. For any admissible number $\zeta$ there exists a unique central bicharacter $\bc_\zeta$ of $\Oo(G_q)$ such that
\begin{equation}\label{3.2}
\bc_\zeta (u^i_j\otimes u^n_m)=\delta_{ij}\delta_{nm}\zeta,\quad i,j,n,m=1,{\dots},N.
\end{equation}
In order to prove this, one first extends (\ref{3.2}) to a linear functional 
on $\bbbc\langle u^i_j\rangle \otimes \bbbc\langle u^i_j\rangle$ such 
that (CB.1) is satisfied, where $\bbbc\langle u^i_j\rangle$ denotes 
the free bialgebra with generators $u^i_j, i,j=1,{\dots},N$. Because 
$\zeta$ is assumed to be admissible, it follows from the defining 
relations for the algebra $\Oo(G_q)$ that this functional passes to a 
functional $\bc_\zeta$ of $\Oo(G_q)\otimes \Oo(G_q)$. It is easily seen 
that $\bc_\zeta$ is a central bicharacter of $\Oo(G_q)$.

\mn
We now fix a universal $r$-form $\br_{z_0}$ of $\Oo(G_q)$ with parameter $z_0$ as above and denote it by $\br$. Then it is clear that any universal $r$-form $\br_z$ of $\Oo(G_q)$ given by (\ref{3.1}) is of the form $\br_z=\bc_\zeta\ast\br$ for some admissible number $\zeta$.

\mn
{\bf Proposition 3.4} {\it Suppose that $q$ is not a root of unity. For any universal $r$-form $\bs$ of $\Oo(G_q)$ there exists an admissible complex number $\zeta$ such that $\bs=\bc_\zeta\ast\br$ or $\bs=\bc_\zeta\ast\bar{\br}_{21}$. The coquasitriangular Hopf algebra $\Oo(G_q)$ has an essentially unique universal $r$-form.}

\mn
{\bf Proof.} Suppose that $\bs$  is a universal $r$-form of $\Oo(G_q)$. We define $N^2\times N^2$-matrices $T=(T^{in}_{jm})$ and $\hat{T}=(\hat{T}^{ni}_{jm})$ by
$$
\hat{T}^{ni}_{jm}=T^{in}_{jm}:=s(u^i_j\otimes u^n_m),\quad i,j,n,m=1,{\dots},N.
$$
From the general theory of coquasitriangular Hopf algebras 
(see, for instance, \cite{Klimyk}, 10.1) we conclude that 
$\hat{T}$ belongs to the centralizer algebra $\Mor (u\otimes u)$ of 
the tensor product corepresentation $u\otimes u$ and that $T$ 
satisfies the quantum Yang-Baxter equation, so $\hat{T}$ fulfills the braid relation 
\begin{equation}\label{3.3}
\hat{T}_{12}\hat{T}_{23}\hat{T}_{12}=\hat{T}_{23}\hat{T}_{12}\hat{T}_{23}
\end{equation}
on the space of the tensor product corepresentation $u\otimes u\otimes u$. 

\mn
We shall carry out the proof in the cases $O_q(N)$ and $Sp_q(N)$. Then it 
is well-known (see \cite{Reshetikhin} or \cite{Klimyk},  Proposition 8.40) 
that there exists a homomorphism $\pi$ of the Birman-Wenzl-Murakami 
algebra $\BWM_3(q,\epsilon q^{N-\epsilon})$ to the centralizer algebra 
$\Mor (u\otimes u \otimes u)$ such that $\pi(G_i)=\hat{R}_{i,i+1}, i=1,2$, 
where $\epsilon=1$ for $O_q(N)$ and $\epsilon=-1$ for $Sp_q(N)$. The 
generators of the BWM-algebra $\BWM_3(q,\epsilon q^{N-\epsilon})$ are 
denoted by $G_1, G_2,E_1,E_2$ and their images under the homomorphism 
$\pi$ by $g_1,g_2,e_1,e_2$. The matrices $\{\hat{R}, \hat{R}^{-1}, I\}$ 
form a basis of the vector space $\Mor (u\otimes u)$. Since $\hat{T}\in\Mor(u\otimes u),$ there are complex numbers $\alpha, \beta,\gamma$ such that $\hat{T}=\alpha\hat{R}+\beta\hat{R}^{-1}+{\gamma\cdot I}$. Thus, we have
\begin{equation}\label{3.4}
\hat{T}_{i,i+1}=\alpha g_i+\beta g_i^{-1}+{\gamma\cdot I}.
\end{equation}
The crucial step of the proof is to show that the braid relation for $\hat{T}$ implies that $\hat{T}$ is a complex multiple of either $\hat{R},\hat{R}^{-1}$ or $I$. In order to prove this, we essentially use the relations of the BWM-algebra $\BWM_3(q,\epsilon q^{N-\epsilon})$ (see \cite{Birman}, p. 225).\\
Inserting (\ref{3.4}) into (\ref{3.3}), both sides of (\ref{3.3}) are sums of
27 summands. From the braid  relation $G_1G_2G_1=G_2G_1G_2$ in the BWM-algebra
we get $g_1g_2g_1=g_2g_1g_2$,
$g_1^{-1}g_2^{-1}g_1^{-1}=g_2^{-1}g_1^{-1}g_2^{-1}$,
$g_1^{-1}g_2g_1=g_2g_1g_2^{-1}$, $g_1g_2^{-1}g_1^{1}=g_2^{-1}g_1^{-1}g_2$,
$g_1^{-1}g_2^{-1}g_1=g_2g_1^{-1}g_2^{-1}$ and $g_1g_2g_1^{-1}=g_2^{-1}g_1g_2$.
Inserting these relations and cancelling equal terms on both sides of (\ref{3.3}) we finally obtain

\begin{gather}
\alpha^2 \beta g_1 g^{-1}_2 g_1 + \alpha \beta^2 g_1^{-1} g_2 g^{-1}_1 + 
\alpha^2 \gamma g^2_1 + \beta^2 \gamma g^{-2}_1\notag\\
\label{3.5} =\alpha^2 \beta g_2 g^{-1}_1 g_2 + \alpha \beta^2 g^{-1}_2 g_1 g^{-1}_2 + \alpha^2 \gamma g^2_2 + \beta^2 \gamma g^{-2}_2~.
\end{gather}
Now we recall the following relations in the BWM-algebra 
(see \cite{Birman}, p. 255):
\begin{align*}
&G^{-1}_i-G_i=\lambda E_i-\lambda \cdot 1, G_2 E_1 G_2=G^{-1}_1 E_2 G^{-1}_1, G^{-1}_2 E_1 G^{-1}_2=G_1 E_2 G_1,\\
&E_1E_2G_1=E_1G^{-1}_2, G_1E_2E_1=G^{-1}_2E_1, E_1E_2E_1=E_1,
\end{align*}
where $\lambda:=q-q^{-1}.$ Applying the images of these relations 
under the homomorphism $\pi$, a  straightforward computation shows 
that (\ref{3.5}) reduces to the equation 
\begin{gather}
\alpha^2(\gamma-\beta\lambda)(g^2_1-g^2_2)+
\beta^2(\alpha\lambda+\gamma)(g^{-2}_1-g^{-2}_2)\notag\\
+\alpha\beta(\alpha+\beta)\lambda^2\{e_2g_1+g_1e_2-e_1g^{-1}_2-g^{-1}_2e_1+\lambda(e_1e_2+e_2e_1+e_1+e_2)\}=0.
\label{3.6}
\end{gather}
The BWM-algebra $\BWM_3(q,\epsilon q^{N-\epsilon})$ is 15-dimensional and the elements
\begin{equation}\label{3.7}
1,G_1,G_2,G_1G_2,G_2G_1,G_1G_2G_1,E_1,E_2,E_1E_2,
E_2E_1,G_1E_2,E_2G_1,G^{-1}_2E_1, E_1G^{-1}_2,G_1E_2G_1
\end{equation}
form a vector space basis (see \cite{Birman} or \cite{Klimyk}), Proposition
8.39). From the representation theory of quantum groups it is well-known (see,
for instance, \cite{Klimyk}, 8.6.2) how to decompose the tensor product
$u\otimes u\otimes u$ into irreducible components. From these decompositions
it follows that the homomorphism $\pi$ is injective for $O_q(N),N\ge 3$, and
for $Sp_q(N),N\ge 6$. In these cases the images of the elements (\ref{3.7})
under $\pi$ are also linearly independent. Therefore, the coefficient of,
say, $e_2g_1$ in (\ref{3.6}) is zero, so that 
\begin{equation}\label{3.8}
\alpha\beta(\alpha+\beta)=0.
\end{equation}
Hence the second line of (\ref{3.6}) vanishes identically and the coefficients of $g_1$ and $e_1$ in (\ref{3.6}) have to be zero as well. Using the relation $G_1E_1=\epsilon q^{\epsilon-N}E_1$ in the algebra $\BWM_3(q,\epsilon q^{N-\epsilon})$, these coefficients are computed as
\begin{align}\label{3.9}
&-\alpha^2(\gamma\beta\lambda)\lambda+\beta^2(\alpha\lambda+\gamma)\lambda =0,\\
\label{3.10}
&\alpha^2(\gamma-\beta\lambda)\lambda\epsilon q^{\epsilon-N}+\beta^2(\alpha\lambda+\gamma)\lambda(\epsilon q^{N-\epsilon}-\lambda)=0,
\end{align}
respectively. In the case of the quantum group $Sp_q(4)$ the 
corepresentation corresponding to the Young tableaux of a 
column with 3 boxes does not occur in the decomposition of 
tensor product $u\otimes u\otimes u$ (see \cite{Klimyk}, p.289) and 
we have $\dim \Mor (u\otimes u\otimes u)=14$. An explicit computation 
shows that the images of the basis elements (\ref{3.7}) satisfy the 
linear relation
\begin{align}
&1-q^{-1}g_1-q^{-1} g_2+g^{-2} g_1g_2+q^{-2} g_2g_1-q^{-3} g_1g_2g_1-q^{-6} e_1-q^{-2} e_1-q^{-4} e_1 e_2\notag\\
\label{3.11}
&-q^{-4} e_2 e_1+q^{-3} g_1e_2+q^{-3} e_2g_1+g^{-5}e_1g^{-1}_2
+q^{-5} g^{-1}_2 e_1-q^{-4} g_1 e_2 g_1=0.
\end{align}
Therefore, the derivation of equations ({\ref{3.8})--({\ref{3.10}) from ({\ref{3.6}) is also valid in the case $Sp_q(4)$. Since $(q-\epsilon q^{\epsilon-N})(q+\epsilon q^{\epsilon-N})\ne 0$ by assumption, ({\ref{3.9}) and ({\ref{3.10}) imply that
\begin{equation}\label{3.12}
\alpha^2 (\gamma-\beta\lambda)=\beta^2 (\alpha\lambda+\gamma)=0~.
\end{equation}
The solutions of equations ({\ref{3.8}) and ({\ref{3.12}) are $\beta=\gamma=0, \alpha=\gamma=0$ and $\alpha=\beta=0$ which means that $\hat{T}$ is a multiple of either $\hat{R}, \hat{R}^{-1}$ or $I$.

\sn
For the quantum groups $GL_q(N)$ and $SL_q(N)$ the proof is similar 
and much simpler. Then there is a homomorphism of the Hecke 
algebra $H_3(q)$ on the centralizer algebra 
$\Mor (u\otimes u\otimes u)$ which is injective for $N\ge 3$. 
In the case $N{=}2$ we have $\dim H_3(q)=1+\dim\Mor 
(u\otimes u\otimes u){=}6$ and the corresponding linear relation is 
obtained from ({\ref{3.11}) by setting $e_1=e_2=0$ therein. Since the 
Hopf algebras $\Oo(Sp_q(2))$ and $\Oo(SL_{q^2}(2))$ are isomorphic, we 
also cover the case $Sp_q(2)$ in this manner which was excluded during 
preceding considerations.

\sn
Summarizing, we have shown that $\hat{T}=z\cdot \hat{R}$ or 
$\hat{T}=z\cdot \hat{R}^{-1}$ or $\hat{T}=z\cdot I$ for some complex 
number $z$. The rest of the proof is more or less routine. Using the 
fact that $\bs$ is a dual pairing of $\Oo(G_q)$ and $\Oo(G_q)^{\rm op}$ 
it follows that the case $\hat{T}=z\cdot I$ is impossible (because it is 
not compatible with the defining relation $\hat{R}_1u_1u_2=u_1u_2\hat{R}$) 
and that the number $z$ must be as described at the beginning of this 
section. Thus we have $\bs=\br_z$ or $\bs=(\bar{\br}_z)_{21}$. Fixing 
a universal $r$-form $\br_{z_0}$ and reformulating the latter in terms 
of a central bicharacter $\bc_\zeta$, the proof will be completed. \hfill\qed

\mn
{\bf Remark  3.5} As I have learned from the referee, the universal r-forms of
the bialgebra $\Oo(M_q(N))$ have been described recently in the paper 
\cite{Takeuchi2}. This result implies the assertion of Proposition 3.4 
in the case $G_q=GL_q(N)$.

\vfill\eject\noindent
\bn
{\bf 4. On functionals describing the square of the antipode}

\mn
Throughout this section we assume that $\A$ is a coquasitriangular Hopf algebra and $\br$ is a universal $r$-form of $\A$.\\
Let $f_\br$ and $\bar{f_\br}$ denote the linear functionals on $\A$ defined by 
\begin{equation}\label{4.1}
f_\br(a)=\br(a_{(1)}\otimes S(a_{(2)}))~{\rm  and }~~~ \bar{f_\br} (a)=\bar{\br}(S(a_{(1)})\otimes a_{(2)}),~ a\in\A.
\end{equation}
Then  it is well-known (see, for instance, \cite{Klimyk}, Proposition 10.3) that  $\bar{f_\br}$ is the convolution inverse of $f_\br$ and that the square of the antipode $S$ is given by
\begin{equation}\label{4.2}
S^2=\bar{f_\br}\ast \id\ast f_\br,~~{\rm that~is},~~ S^2 (a)=\bar{f}_\br (a_{(1)}) a_{(2)} f_\br (a_{(3)}), a\in\A.
\end{equation}

\mn
{\bf Lemma 4.1} {\it The functional $f_\br$ satisfies the equations}
\begin{equation}\label{4.3}
\br_{21}\ast\br \ast (f_\br \circ m_\A) = 
(f_\br \circ m_\A)\ast \br_{21}\ast\br = f_\br\otimes f_\br.
\end{equation}

\mn
{\bf Proof.} Using the properties (CQT.1)--(CQT.3) of the universal $r$-form $\br$ we compute
\begin{align*}
\br_{21}(a_{(1)}\otimes &b_{(1)})\br(a_{(2)}\otimes b_{(2)})f_\br (a_{(3)} b_{(3)})\\
&=\br(b_{(1)}\otimes a_{(1)})\br(a_{(2)}\otimes b_{(2)})\br (a_{(3)} b_{(3)}\otimes S(a_{(4)}b_{(4)}))\\
&=\br(b_{(1)}\otimes a_{(1)})\br(a_{(3)}\otimes b_{(3)})\br (b_{(2)} a_{(2)}\otimes S(b_{(4)})S(a_{(4)}))\\
&=\br(a_{(4)}\otimes b_{(4)})\br(b_{(1)}\otimes a_{(1)})\br (b_{(2)} a_{(2)}\otimes S(a_{(5)})))\br(b_{(3)} a_{(3)}\otimes S(b_{(5)})) \\
&=\br(a_{(4)}\otimes b_{(4)})\br(b_{(2)}\otimes a_{(2)})\br (a_{(1)} b_{(1)}\otimes S(a_{(5)})))\br(b_{(3)}\otimes S(b_{(6)}))\br(a_{(3)}\otimes S(b_{(5)}))\\
&=\br(b_{(2)}\otimes a_{(2)})\br(a_{(1)}b_{(1)}\otimes S( a_{(3)}))\br(b_{(3)}\otimes S(b_{(4)})) \\
&=\br(b_{(2)}\otimes a_{(2)})\br(a_{(1)}\otimes S(a_{(4)}))\br (b_{(1)} \otimes S(a_{(3)}))f_\br(b_{(3)})\\
&=f_\br(a)f_\br(b)
\end{align*}
for $a,b \in \A$. This proves the first equality 
$\br_{21}\ast\br \ast (f_\br \circ m_\A) = f_\br\otimes f_\br.$ 
Applying condition (CQT.3) twice to this relation we get the second equality
$(f_\br \circ m_\A)\ast \br_{21}\ast\br = f_\br\otimes f_\br.$
\hfill\qed

\sn
In some sense the functional $\br_{21}\ast\br$ on $\A \otimes \A$ 
measures the distance of $f_\br$ from being a character. In particular, Lemma 3.1 implies 

\bn
{\bf Corollary 4.2}. {\it The functional $f_\br$ is a character of $\A$ (that is, $f_\br (ab)=f_\br(a) f_\br(b)$ for $a,b\in \A$) if and only if $\A$ is cotriangular (that is, $\bar{\br}=\br_{21}$).}

\bn
Recall that any coquasitriangular Hopf algebra $\A$ has a second universal $r$-form $\bs :=\bar{\br}_{21}$ given by $\bs (a\otimes b)=\bar{\br} (b\otimes a), a,b\in\A$.

\mn
{\bf Proposition 4.3} {\it The functionals $f_\br,\bar{f_\br}, f_\bs,\bar{f_\bs}$ pairwise commute in the algebra $\A^\circ,z{:=}f_\br\ast \bar{f_\bs}$ belongs 
to the center of $\A^\prime$ and  $g:=f_\br\ast f_\bs$ is a 
character of $\A$ such that $S^4=\bar{g}\ast \id\ast g$.}

\mn
{\bf Proof.} The antipode $S$ is bijective and we have $S\ast \bar{f}_\bs=\bar{f_\bs}\ast S^{-1}$ (\cite{Klimyk}, Proposition 10.3). Using this fact, the relation $\br (S(a)\otimes S(b))=\br(a\otimes b)$ and (\ref{3.2}), we obtain
\begin{align*}
f_\br\ast \bar{f_\bs} (a)&=\br (a_{(1)}\otimes S(a_{(2)}))\bar{f_\bs }(a_{(3)})=\br (a_{(1)}\otimes S^{-1}(a_{(3)}))\bar{f_\bs} (a_{(2)})\\
&=\br (S^2(a_{(1)})\otimes S(a_{(3)}))\bar{f_\bs} (a_{(2)})
=\br (a_{(2)}\otimes S(a_{(3)}))\bar{f_\bs} (a_{(1)})=\bar{f_\bs} \ast f_\br (a).\\
\end{align*}
Hence the functionals $f_\br, \bar{f_\br}, f_\bs, \bar{f_\bs}$ pairwise commute. \\
Applying (\ref{3.2}) to both $\br$ and $\bs$, we get 
$$
z\ast \id=f_{\br}\ast \bar{f_\bs}\ast \id = f_\br \ast S^2\ast\bar{f_\bs}= \id \ast f_\br \ast \bar{f_\bs}=\id \ast z,
$$
so $z$ is in the center of the algebra $\A^\prime$.\\
Since $\br_{21}\ast\br\ast \bs_{21}\ast\bs=\varepsilon_{\A\otimes\A}$, 
the equations (\ref{4.3}) easily imply that $g(ab) = g(a)g(b)$ for $a,b\in\A$, that is, $g$ is a character on $\A$. By (\ref{4.2}), we have $S^4=\bar{g}\ast\id\ast g$.\hfill\qed

\mn
We illustrate the preceding by a simple example.

\bn
{\bf Example 4.4} Let $\A$ be the Hopf algebra $\bbbc\bbbz$ of the  group of integers. Then any universal $r$-form of $\A$ is of the form $\br(n\otimes m)=\lambda^{nm}, n,m\in\bbbz$, for some fixed $\lambda\in\bbbc, \lambda\ne 0$ (see also Example 3.3). Thus we have $f_\br(n)=\lambda^{n^2}, f_\bs(n)=\lambda^{-n^2}$ and $F_\br(n)=1$ for $n\in\bbbz$, so that $F_\br=\varepsilon$. Since $f_\br(n+m)=f_\br(n)f_\br(m)\lambda^{2nm}$ for $n,m\in\bbbz$, the functional $f_\br$ on $\A$ is far from being a character in general.

\mn
Let us suppose now that $\A$ is a  coquasitriangular Hopf algebra equipped with a universal $r$-form $\br$ and that $\A$ is also a $CQG$-algebra (see \cite{DK} or \cite{Klimyk}, 11.3.1, for this notion). Let $f_z$, $z\in\bbbc$, denote Woronowicz$^\prime$s modular characters on $\A$ (see \cite{Woronowicz} or \cite{Klimyk}, 11.3.4). Recall that $f_z\ast f_{z^\prime}=f_{z+z^\prime}$, for $z,z^\prime\in\bbbc$ and $S^2(a)=f_2\ast\id\ast f_{-2}$ for $a\in\A$. Then the functionals $F_\br$ and $f_{-2}$ are both characters on $\A$ which implement $S^4$, that is, 
$$
S^4(a)=\bar{F}_\br(a_{(1)}) a_{(2)} F_\br(a_{(3)})=f_2(a_{(1)})a_{(2)} f_{-2} (a_{(3)}),~a\in\A.
$$
Hence $F_\br f_2$ is a character of $\A$ which is central in $\A^\prime$. This suggests the following

\bn
{\bf PROBLEM:} ~{\it Do the characters $F_\br$ and $f_{-2}$ on $\A$ coincide?}

\mn
If the Hopf algebra $\A$ is cocommutative (in particular, if $\A$ is the group algebra $\bbbc G$ of an abelian group $G$), then we have $F_\br=f_{-2}=\varepsilon$ and so the answer is affirmative.\\
A more interesting case is the Hopf $\ast$-algebra $\Oo(G_q)$, where $G_q$ is
one of the compact forms $U_q(N)$, $SU_q(N)$, $O_q(N,\bbbr)$, $USp_q(N)$
of the quantum groups $GL_q(N)$, $SL_q(N)$, $O_q(N)$, $Sp_q(N)$, respectively,
and $q$ is real. Then $\Oo(G_q)$ is a $CQG$-algebra (\cite{Klimyk},
Example 11.7) and a coquasitriangular Hopf algebra with universal $r$-form
$\br$ determined by (\ref{3.1}). 

\bn
{\bf Proposition 4.5} {\it Then we have $F_\br=f_{-2}$.}

\bn
{\bf Proof.} From the explicit formulas for $f_\br (u^i_j),f_\bs(u^i_j)$ and $f_{-2}(u^i_j)$ listed in \cite{Klimyk}, p.341 and p.425, respectively, we see that $F_{\br}(u^i_j)=f_{-2}(u^i_j)$ for $i,j=1,{\dots}, N$. Therefore, since $F_\br$ and $f_{-2}$ are both characters, they coincide on the whole algebra $\Oo(G_q)$.\hfill\qed

\newpage\noindent
{\sc Acknowledgement:}
I would like to thank Drs. I. Heckenberger and A. Sch\"uler 
and Prof. M. Takeuchi 
for their useful comments on this paper.

\mn

\end{document}